
\input amstex

\magnification=1200
\loadmsam
\loadmsbm
\loadeufm
\loadeusm
\UseAMSsymbols

\hsize=6.0truein
\hoffset=0.15truein
\vsize=9truein
\voffset=-0.2truein

\def\leftitem#1{\item{\hbox to\parindent{\enspace#1\hfill}}}

\def\boxit#1#2{\hbox{\vrule
	\vtop{%
	\vbox{\hrule\kern#1%
	\hbox{\kern#1#2\kern#1}}%
	\kern#1\hrule}%
	\vrule}}

\def\leaderfill{\leaders\hbox to 1em{\hss.\hss}\hfill}

\parskip=\medskipamount
\document

\input epsf

\centerline{\bf Torsion Elements in  the Mapping Class Group of a Surface}

\bigskip

\centerline{\bf Feng Luo}

\bigskip

{\bf Abstract}
Given a finite set of $r$ points in a closed surface of genus $g$, 
we consider the torsion elements in the mapping class group of the surface 
leaving the finite set invariant. We show that the torsion elements generate
the mapping class group if and only if  $(g, r) \neq (2, 5k+4)$ for
some integer $k$.
\bigskip

\S1. {\bf Introduction}

1.1.
The purpose of this paper is to investigate when the mapping class group
of a compact surface is generated by torsion elements.  Our result
gives a complete answer to this question. 

{\bf Theorem.} \it Suppose $\Sigma_{g,r}$ is a compact orientable
surface of genus $g$ with $r$ many boundary components where $g,r
\geq 0$. Then the mapping class group of the path
components of the orientation preserving homeomorphisms of
the surface  is generated
by torsion elements if and only if $(g, r) \neq (2, 5k+4)$ for some
integer $k \in \bold Z$. The torsion elements in the mapping class group of 
the surface $\Sigma_{2, 5k+4}$ generate an index 5 subgroup.

Furthermore, in the case of $(g, r) \neq (2, 5k+4)$, the order $n$ of
the torsion elements generating the group can be chosen as follows:
(a) if $g \geq 3$, then $n=2$; (b) if $g=2$, then $n \in \{2,5\}$;
(c) if $g=1$, then $n \in \{2,3, 4,\}$; and (d) if $g=0$, then
$n \in \{r-1, r\}$. \rm

Note that by identifying each boundary component to a point, one sees
that the mapping class group in the theorem is the same as the
mapping class group of  a closed surface leaving a set of $r$ points
invariant.

The existence of the exceptional cases $(2, 5k+4)$ is caused by the fact
that there is only one non-trivial $\bold Z_5 $-action on the
closed surface of genus-2. And this $\bold Z_5 $-action has
too few fixed points (only 3 points). 

The result is motivated by the example
of the torus whose mapping class group is $SL(2, \bold Z) = 
\bold Z_4 *_{\bold Z_2} \bold Z_6$.
In the case of surfaces of genus at least 3, the theorem can be derived
easily from the work of Harer [Ha]. A simple derivation of it will be
given in \S1.3. The main body of the paper is to prove the theorem for surfaces
of genus at most 2.

On the related question of Torelli groups, the following can be derived
easily from the work of Johnson [Jo] and Powell [Po]. Namely the
Torelli group of a closed surface of genus at least 3 is generated by
the products of even number of hyperelliptic involutions.

1.2. 
A basic strategy to prove that a group $G$ is generated by torsion
is to produce a short exact sequence $1 \to H \to G \to G/H \to 1$ 
so that $H$ is generated by torsion elements in $G$ and there are torsion elements
in $G$ whose projections in $G/H$ generate the group. In our case,
let $\Gamma^*_{g,r}$ be the mapping class group of the surface
which is the group of path
components of orientation preserving homeomorphisms of the surface.
We use the short exact sequence $ 1 \to [\Gamma^*_{g,r}, \Gamma^*_{g,r}]
\to \Gamma^*_{g,r} \to H_1(\Gamma^*_{g,r}) \to 1$. We first
show that the mapping class group is generated by
Dehn-twists and torsion elements. Then we use it to prove
that the commutator subgroup $[\Gamma^*_{g,r}, \Gamma^*_{g,r}]$
is generated by torsion elements. Finally we show that
the images of the torsion elements in the first homology
group generate the homology.

The most interesting  part of the proof is to show that for genus-2
surface $\Sigma_{2,0}$ the projections of 
a 5-fold symmetry and an involution with two fixed points in the
surface generate the first homology group
$H_1(\Gamma^*_{2,0} ) \simeq$ $\bold Z_{10}$.

In the process of proof, we show that
each compact surface admits a non-trivial $\bold Z_3 $ action.
We also calculated the first homology of $\Gamma^*_{g,r}$ using
Harer and Powell's work.

1.3. Using the work of Harer [Ha], we give a simple proof 
the theorem for surfaces of genus at least 3 in this section. 
The argument is no longer working for low genus
surfaces

Let $\Gamma_{g,r}$ be the \it pure mapping class group \rm
$\Gamma_{g,r}$ which is the subgroup of $\Gamma^*_{g,r}$ so that it
acts trivially on  the set of boundary components of the surface.
Then there is a short exact sequence
$1 \to \Gamma_{g,r} \to \Gamma^*_{g,r} \to S_r \to 1$ where $S_r$
is the permutation group on $r$ boundary components of the surface.
It is easy to show that
there are torsion elements in the mapping class group $\Gamma^*_{g,r}$
whose projections in $S_r$ generate the permutation group. Thus
it remains to show that the pure mapping class group $\Gamma_{g,r}$ 
is generated by torsion elements in $\Gamma^*_{g,r}$.  Note that
if $r > 2g+2$, then the pure mapping class group $\Gamma_{g,r}$
is torsion free.

The fundamental work of Dehn [De] and Lickorish [Li] shows that 
$\Gamma_{g,r}$ is generated by Dehn-twists. If the genus of
the surface is positive, it can be shown that Dehn-twists on
non-separating simple loops generate the pure mapping class group.
On the other hand, any two Dehn-twists on non-separating simple
loops are conjugate. Thus it suffices to show that the Dehn-twist
on a non-separating simple loop is a product of torsions in
$\Gamma^*_{g,r}$. 

Since the genus of the surface is at least three, there is
a 4-holed sphere subsurface in $\Sigma_{g,r}$ bounded by
four non-separating simple loops $a,b,c,d$ so that the
complement of the 4-holed sphere is connected. Take
three simple loops $x,y,$ and $z$ in the 4-holed sphere
forming a lantern position (see figure 4.1(a)). Then the lantern relation gives:
$ABCD = XYZ$ where capital letters are the positive Dehn-twist
on the small letter simple loops. Thus $A = XB^{-1} YC^{-1} ZD^{-1}$
where each of $(x,b)$, $(y,c)$ and $(z,d)$ is a pair of
disjoint non-separating simple loops so that complement of
their union in $\Sigma_{g,r}$ is connected. Thus it suffices to
show, for instance, that $XB^{-1}$ is a product of involutions.
By the choice of $(x,b)$, there is an involution $f$ of the surface
sending $x$ to $b$. Thus $B = fXf^{-1}$. This shows that
$XB^{-1} = (X f X^{-1}) f^{-1}$ which is the product of two
involutions $f^{-1}$ and $X f X^{-1}$. In Harer's proof [Ha],
he used the equation $XB^{-1} = X f X^{-1}f^{-1}$ to show
that $A$ is a product of commutators  where $f$ is an element
in $\Gamma_{g,r}$ and deduced that the first
homology of $\Gamma_{g,r}$ is trivial for $g \geq 3$.

The proof for the Torelli group of closed surface of genus at least 3
is the same as above using the work of Johnson [Jo] and Powell [Po].
Indeed, they showed that the Torelli group is generated by a product
of Dehn-twists $AB^{-1}$ where the simple loops $a,b$ are disjoint, both
non-separating and $a \cup b$ decomposes the surface into two pieces.
Thus there is a hyperelliptic involution $f$ sending $a$ to $b$ and
$AB^{-1} = (AfA^{-1})f^{-1}$.

1.4.  Periodic homeomorphisms and Dehn-twists are the
two extremal cases of self homeomorphisms in the sense that
the Dehn-twists have the largest support and periodic homeomorphisms
have the smallest support. In particular this makes the
calculation of a periodic homeomorphism in terms of Dehn-twists
quite complicated.
In view of the fact that
the mapping class group $\Gamma_{1,0}$ of a torus is
$SL(2, \bold Z) =
\bold Z_4 *_{\bold Z_2} \bold Z_6$ and also the recent work 
of Wajnryb  [Wa] that the pure mapping class group is generated
by two elements, one is attempting  to ask what is the minimal
number of periodic generators for the mapping class group. Is it
possible to generate $\Gamma^*_{g,r}$ by torsion elements where 
the number of generators is independent of  $g,r$? 

1.5.  The organization of the paper is as follows. In section 2, we 
recall some well known symmetries of surfaces of low genus and
their expressions in Dehn-twists. In section 3 we show that
the commutator subgroup of the mapping class group is generated by
torsion elements. In section 4, we prove the theorem using
the first homology group of the mapping class group.

1.6. Acknowledgment. The work is supported in part by the NSF.

\S2. {\bf Period Homeomorphisms on Low Genus Surfaces}

We introduce notations and terminologies in this section.
Also we recall some well known symmetries of surfaces of genus at most
2.

2.1. We start by introducing some notations and conventions.
Given a finite set $X$, we use $|X|$ to denote the number of
elements in $X$.
Surfaces are assumed to be oriented. Subsurfaces have the induced
orientation. Simple loops on surfaces will be denoted  by $a$, $b$, ...,$c$.
Positive Dehn-twists on them will be denoted by $D_a$, $D_b$, ..., $D_c$
or simply by $A$, $B$, ..., $C$ if no confusion will arise. 
A small regular neighborhood of a 
1-dimensional submanifold $s$ will be denoted by $N(s)$. 
We shall not distinguish homeomorphisms from their isotopy classes.
Thus a Dehn-twist sometimes means the isotopy class of a Dehn-twist.


Given an arc $s$ joining two boundary components $\partial_1$ and
$\partial_2$ of a surface, the  \it positive half-Dehn-twist \rm
(or simply half-twist)
along $s$ will be denoted by $D^{1/2}_s$. It is a self-homeomorphism
$h$ supported in $N(s \cup \partial_1 \cup \partial_2)$ so that
$h$ leaves $s$ invariant and interchanges $\partial_1$, $\partial_2$
and $h^2$ is the positive Dehn-twist along $\partial N(s \cup
\partial_1 \cup \partial_2) -\partial_1 \cup \partial_2$.

Given a compact surface $\Sigma$, let $\Sigma^*$ be the quotient of
$\Sigma$ obtained by identifying each boundary component to a point.
Clearly each homeomorphism $h$ of $\Sigma$ induces a canonical
homeomorphism $h^*$ of the closed surface $\Sigma^*$.
An 
\it hyperelliptic involution \rm $h$ on a surface $\Sigma$ is
an involution so that its induced involution $h^*$ has exactly $2g+2$ 
fixed points in $\Sigma^*$ where $g$ is the genus. Note that
any two hyperelliptic involutions on a closed surface are conjugate.
Also the action of a hyperelliptic involution on the first homology
of a closed surface is the multiplication by $-1$.

2.2. The fundamental work of 
Dehn [De] and Lickorich [Li] states that the pure mapping
class group  $\Gamma_{g,r}$ is finitely generated by Dehn-twists. 
We shall need a
slightly improved version of it. See [Ge] for a proof.

{\bf Proposition.} \it ([Ge]) For a surface of positive genus, the pure mapping
class group is finitely generated by Dehn-twists on non-separating
simple loops $a_1, ..., a_n$ so that for each pair of indices $i,j$, either
$a_i$ is disjoint from $a_j$ or $a_i$ intersects $a_j$ at one point.
\rm

For instance the following loops form such a generating set. See figure 2.1.

\midspace{0.1cm}
\centerline{\epsfbox{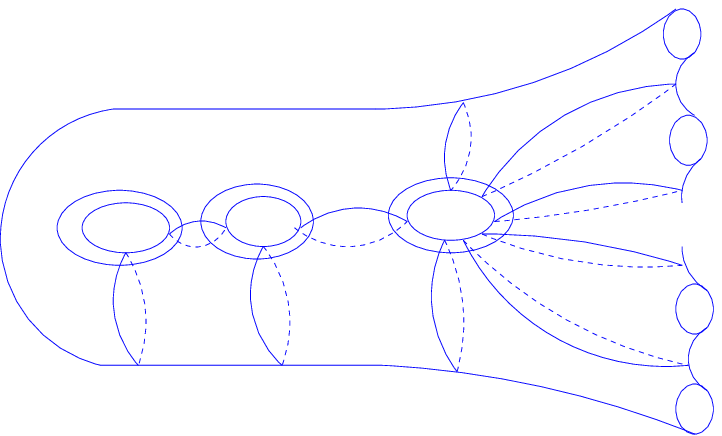}}
\midspace{0.1cm}
\centerline{Figure 2.1}

2.3. The basic symmetries of surfaces of low genus that we shall use are
the followings. The $2\pi/n$-rotation $\tau_n$ of the 2-sphere,
the $\bold Z_4$ and $\bold Z_6$ non-free actions on the torus and
the $\bold Z_5$ action on the genus-2 surface.

The action of $\bold Z_4$ on the torus is constructed as follows.
Consider the torus as the quotient of the square by identifying the
opposite sides. Then the generator $\tau_4$ of the $\bold Z_4$
action is induced by the $\pi/2$-rotation of the square.
It has 2 fixed points corresponding to the center and vertices and
has an orbit consisting of two points coming from the
mid-points of the four sides. Note that $\tau_4^2 = \tau_2$ is a
hyperelliptic involution.

The action of $\bold Z_6$ on the torus is constructed as follows.
Consider the torus as the quotient of the hexagon by identifying the
opposite sides. Then the generator $\tau_6$ of the action is
induced by the $2\pi/6$-rotation of the hexagon.
It has  one fixed point corresponding to the center of the hexagon.
And it has an orbit consisting of three points coming from the
mid-points of the six sides and an orbit consisting of two points
coming from the  six vertices. Note that $\tau_6^3 = \tau_2$ is a
hyperelliptic involution and that $\tau_6^2 =\tau_3$ is
a $\bold Z_3$ action having three fixed points.

The action of $\bold Z_5$ on the genus-2 surface is induced by the
$2 \pi/5$-rotation of the 10-sided polygon. Here the genus-2 surface
is the quotient of the 10-sided polygon by identifying the opposite
sides (in fact this shows that $\bold Z_{10}$ acts on the surface).
The generator $\tau_5$ of the action has three fixed points in
the genus-2 surface. These fixed points are the quotients of
the center and the vertices of the 10-sided polygon.

{\bf Lemma.} \it (a) For $k=2,3,4,6$, each genus-1 surface $\Sigma_{1,r}$ 
has an order $k$ periodic  homeomorphism $f_k$  so that
 its induced homeomorphism in the torus $\Sigma_{1,0} =\Sigma_{1,r}^*$ is
the standard symmetry $\tau_k$.

(b) For $r \neq 5k+4$, the surface $\Sigma_{2, r}$ has an order
$5$ periodic homeomorphism $f_5$ so that its induced homeomorphism
on the genus-2 surface is the standard symmetry $\tau_5$.
There is no $\bold Z_5$ action on the surface $\Sigma_{2, 5k+4}$.

(c) Suppose the genus $g$ is positive. If $0 \leq k \leq  3$  and
$r-k$ is a non-negative even integer,
then
there is an involution on the surface $\Sigma_{g,r}$ so that exactly $k$ 
many boundary components of the surface is invariant under the involution.

(d) If $\tau$ is an involution on a closed surface $\Sigma_{g,0}$ with at 
least one fixed point, then there exists an involution in each
compact surface $\Sigma_{g,r}$ which induces $\tau$ in the
closed surface $\Sigma_{g,0} = \Sigma_{g,r}^*$.

\rm

{\bf Proof.} A simple method to construct a period-$n$ homeomorphism
$f_n$ on the surface $\Sigma_{g,r}$ inducing a given period-$n$
homeomorphism $h_n$ on the closed surface $\Sigma_{g,0} = \Sigma_{g,r}^*$
is as follows. We take a union of $h_n$-orbits consisting of $r$ points in
$\Sigma_{g,0}$.
Let $\Sigma_{g,r}$ be the surface  obtained from $\Sigma_{g,0}$ by
removing an equivariant  neighborhood of the union of the orbits and
take $f_n$ to be the restriction.  Thus the main issue  is a simple
arithmetic problem.

To see part (a), since $\tau_2$ and $\tau_3$ are powers of $\tau_4$ and
$\tau_6$, it suffices to show the result for $\tau_4$ and $\tau_6$.

For $\tau_4$ on the torus, the homeomorphism has two fixed points $\{a_1, a_2\}$
and an orbit $\{b_1, b_2\}$ consisting of two points. We choose
the union $X$ of $\tau_4$ orbit consisting of $r$ many points
as follows.
If $r \equiv 0$ mod $4$, then $X \cap \{a_1, a_2, b_1, b_2\} = \emptyset$.
If $r \equiv 1,2,3$ mod $4$, then $X \cap \{a_1, a_2, b_1, b_2\}$
$= \{a_1\}$, $\{a_1, a_2\}$ and $\{a_1, b_1, b_2\}$ respectively.

For $\tau_6$ on the torus, the homeomorphism has one fixed point
$a$, an orbit $\{b_1, b_2\}$ consisting of two points, and
an orbit $\{c_1, c_2, c_3\}$ consisting of three points.
We choose the union $X$ of the orbit consisting of $r$ many
points as follows. If $r \equiv 0$ mod $6$, $X \cap \{a, b_1, b_2, c_1,
c_2, c_3\}$ is empty. If $r \equiv 1,2,3,4,5$ mod $6$, we choose
$X$ so that $X \cap \{a, b_1, b_2, c_1, c_2, c_3\}$ is
$\{a\}$, $\{b_1, b_2\}$, $\{c_1, c_2, c_3\}$, $\{a, c_1, c_2, c_3\}$
and $\{b_1, b_2, c_1, c_2, c_3\}$ respectively.

To see part (b),  we note that $\tau_5$ has three fixed points
$\{a_1, a_2, a_3\}$ in $\Sigma_{2,0}$. 
Thus, by the same argument as above, we can construct period-5
homeomorphisms on $\Sigma_{2, 5k+s}$ when $ 0 \leq s \leq 3$. But
there are no $\bold Z_5$ action on $\Sigma_{2, 5k+4}$. Indeed, if 
$f_5$ is such a period-5 homeomorphism, then at least four
boundary components of $\Sigma_{2, 5k+4}$ are invariant under
$f_5$. Thus the induced period-5 homeomorphism $f^*_5$ on the
genus-2 surface has $n \geq 4$ fixed points.
By removing
an equivariant neighborhood of these fixed points, we obtain
a free $\bold Z_5$ action on the surface $\Sigma_{2, n}$ whose
quotient is a surface $\Sigma_{g,n}$ having $n$ boundary components.
But the Euler number multiplies under covering,
i.e., $ -2-n = 5(2-2g-n)$ where $n \geq 4$.
A simple inspection  shows that this is impossible.

Using the same method, one proves part (d).

To prove part (c), the hyperelliptic involution $\tau_2$ on $\Sigma_{g,0}$
has $2g+2 \geq 4$ fixed points. Since $r-k \geq 0$ is even and $k \leq 3$,
we can take  a union $X$ of $\tau_3$ orbits so that
$X$ contains exactly $k$ fixed points of $\tau_2$ and $|X| = r$. 
Now take
$\Sigma_{g,r}$ to be $\Sigma_{g,0}$ with a $\tau_2$ equivariant
neighborhood of $X$ removed. The restriction of $\tau$ gives the
required involution.
$\square$

2.4. In this section, we will express the standard symmetries
in \S2.3 in terms of Dehn-twists and half-twists. 

{\bf Lemma.} \it
(a) ([De]) Let $a$ and $b$ be two simple loops in the torus $\Sigma_{1,0}$
so that they intersect transversely at one point. Let $A$ and
$B$ be the positive Dehn-twists on $a$ and $b$ respectively.
Then the standard symmetries of the torus are the following:
the hyperelliptic involution $\tau_2 = ABABAB$,
the 4-fold symmetry $\tau_4 = ABA$ and the 6-fold symmetry
$\tau_6 = AB$.

(b) (see [Bi]) Let $a_1$, ..., $a_{r-1}$  be the pairwise disjoint
arcs in the planar surface $\Sigma_{0,r}$ so that $a_i$ joins the 
$i$-th boundary $B_i$ to $B_{i+1}$. Let $A_i$ be the half-twist
about the arc $a_i$.
Then $\tau_r = A_1...A_{r-1}$ and $\tau_{r-1} = A_1...A_{r-2}$ 
are $2 \pi/r$ and $2 \pi/(r-1)$-rotation of the surface sending
$a_i$ to $a_{i+1}$ for $1 \leq i \leq r-3$.

(c) (see [Bi]) Let $C_1$,..., $C_5$ be the positive Dehn-twists
on the five simple loops $c_1$, ..., $c_5$ in the genus-2 surface
(see figure 2.3).
Then the hyperelliptic involution $\tau_2 = C_1C_2C_3C_4C_5^2
C_4C_3C_2C_1$ and the 5-fold symmetry is $\tau_5  =\tau_2
C_1C_2C_3C_4$.  There is an involution $\delta$ with
 two fixed points on $\Sigma_{2,0}$ which is a product of
 15 or 25 many Dehn-twists on non-separating simple loops.

\rm

{\bf Proof.}
To see part (a), by the work of Dehn [De],  we have $ABA = BAB$
and also $(ABA)^4 =1$  and $(AB)^6 =1$
(this can be verified using the first homology
group of the torus). Thus the result follows.

To see part (b), we note that a presentation of the mapping
class group $\Gamma^*_{0,r}$ was obtained in [Bi] where the
generators are $A_1, ..., A_{r-1}$ subject to the relations
(1) $ A_i A_j = A_i A_j$ if $|i-j| \geq 2$;
(2) $ A_i A_{i+1} A_i = A_{i+1}A_i A_{i+1}$;
(3) $A_1 ... A_{r-2} A_{r-1}^2 A_{r-2} ... A_1 =1$;
and (4) $(A_1 ... A_{r-1})^n =1$.
The last relation says that $\tau_r = A_1...A_{r-1}$ is a
$2 \pi/r$-rotation on the planar surface. One can also see this using
braid representation of the mapping class group where the half-twist
$A_i$ corresponds to the standard $i$-th  switching generator
$\sigma_i$ of the braid group on $r$-strings. Now the
product $\sigma_1 \sigma_2 ...\sigma_{r-1}$ represents a full twist
string as shown in figure 2.2. Thus $\tau_r =
A_1...A_{r-1}$ corresponds to the $2 \pi/r$-rotation. This also shows
that the product $\tau_j = A_1 A_2 ... A_{j-1}$ leaves the
subsurface $ N(B_1 \cup ... \cup B_{j} \cup a_1 \cup ... \cup a_{j-1})$
invariant and its restriction to the subsurface is a $2 \pi/j$-rotation.  
Now take $j =r-1$ and observe that
the regular neighborhood of  $B_1 \cup B_2 \cup ... \cup B_{r-1} \cup a_1 \cup ... \cup a_{r-2}$ is isotopic to $\Sigma_{0,r}$. We conclude that
$\tau_{r-1}$ is a $2 \pi/(r-1)$-rotation on the surface $\Sigma_{0,r}$.

\midspace{0.1cm}
\centerline{\epsfbox{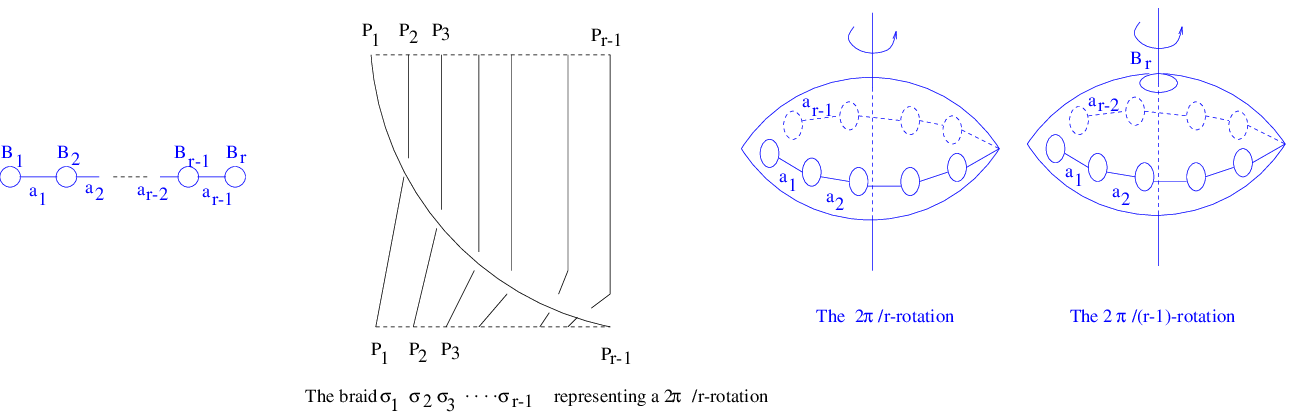}}
\midspace{0.1cm}
\centerline{Figure 2.2}

To see part (c) we recall that the work of
[Li] shows that the mapping class group of the genus-2 closed
surface is generated by five Dehn-twists $C_1, ..., C_5$ on
five simple loops $c_1, ..., c_5$ (see figure 2.3). But these
five simple loops are invariant under the hyperelliptic involution
$\tau_2$ induced by the $\pi$-rotation about the x-axis. Thus
The hyperelliptic
involution $\tau_2$ is in the center of the mapping class group
$\Gamma_{2,0}$ ([Bi], [Vi]). The work of [Bi] and [BH] shows
that the quotient of $\Gamma_{2,0}$ by $\tau_2$ is
the mapping class group $\Gamma^*_{0,6}$. Thus there is a central
extension $ 1 \to \bold Z_2 \to \Gamma_{2,0} \to \Gamma^*_{0,6} \to 1$.
The quotient homomorphism is constructed as follows. Take an orientation
preserving homeomorphism $f$ on $\Sigma_{2,0}$. One may isotope $f$ so
that $f( \tau_2(x)) = \tau_2(f(x)) $ for all $x$.
Thus $f$ induces a homeomorphism $f_*$ on the quotient space
$\Sigma_{2,0}/\tau_2$ which is a 2-sphere with six cone points.
Thus one may think of $f_*$ as an element in $\Gamma^*_{0,6}$. The work
of [Bi] and [BH] shows that the map sending $[f]$ to $[f_*]$ is a well defined
epimorphism with kernel generated by the hyperelliptic involution.
In particular, we see that
the quotient
homomorphism sends the Dehn-twist $C_i$ about the simple loop $c_i$
to the half-twist $A_i$ about $a_i$ in $\Sigma_{0,6}$. 

\midspace{0.1cm}
\centerline{\epsfbox{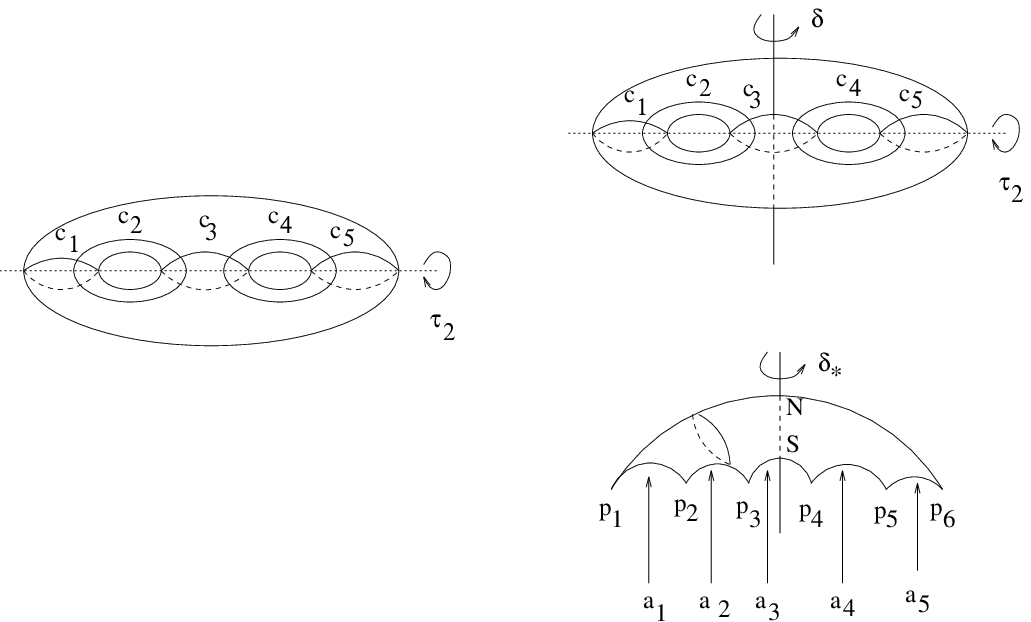}}
\midspace{0.1cm}
\centerline{Figure 2.3}

The lifts
of relation (2)  in $\Gamma^*_{0,6}$ shows that the hyperelliptic
involution $ \tau_2 = C_1C_2C_3C_4C_5^2
C_4C_3C_2C_1$. By part (b), the element $A_1A_2A_3A_4$
 is a period-4  element in $\Gamma^*_{0,6}$. It has two lifts
in $\Gamma_{2,0}$ given by
$C_1C_2C_3C_4$  and $\tau_2 C_1C_2C_3C_4$. The tenth-power of
both lifts are the identity. We claim that $\tau_2
C_1C_2C_3C_4$ has order 5 and is thus equal to $\tau_5$. To this end,
we use the first homology of the surface
$\Sigma_{2,0}$. Note that $H_1(\Sigma_{2,0})$ is generated by
the homology classes $[c_1],..., [c_4]$ where we choose orientation
on $c_i$ so that the algebraic intersection number from $c_i$ to
$c_{i+1}$ is one. A simple calculation shows that  
matrix representative of $C_1C_2C_3C_4$ with respect to the
basis $([c_1],...,[c_4])$ is
$$
\pmatrix
0&1&0&0\\
0&0&1&0\\
0&0&0&1\\
-1&1&-1&1
\endmatrix \right).
$$
The fifth power of the matrix is $-id$. Since the hyperelliptic involution
induces the multiplication by $-1$ in homology, thus the fifth power of
the periodic element $\tau_2 C_1C_2C_3C_4$ is the identity. By
Hurwitz theorem that first homology detects periodic homeomorphisms,
we see that $\tau_5 = \tau_2 C_1C_2C_3C_4$.

To show that there is an involution $\delta$ with two fixed points in
 $\Sigma_{2,0}$ which is
a product of 15 or 25 many Dehn-twists on $c_1, ..., c_5$, we draw the
surface $\Sigma_{2,0}$ as in figure 2.3 where the hyperelliptic involution
$\tau_2$ is induced by the $\pi$-rotation about the x-axis and 
$\delta$ is induced by the $\pi$-rotation about the z-axis.
We claim that $\delta$ is isotopic to a product of 
15 or 25 Dehn-twists on $c_1, ..., c_5$. To see this, let us consider the
induced homeomorphism $\delta_*$ in the quotient sphere $\Sigma_{2,0}/\tau_2$.
The induced involution $\delta_*$ has two fixed points $\{N, S\}$
and leaves
the six cone points $\{p_1,...,p_6\}$ invariant. 
Note that the images of the simple loop
$c_i$'s are the arcs $a_i$ joining $i$-th cone point $p_i$ to $p_{i+1}$.
We construct a homeomorphism from the quotient space $\Sigma_{2,0}/\tau_2$ to the Riemann sphere sending $N$ to the infinity and $S$ to the origin,
the arc $a_1 \cup ... \cup c_5$ into the real axis, and the six
points $(p_1, ..., p_6)$ to $(-3,-2,-1,1,2,3)$ so that the
involution $\delta_*$ becomes $ z \to -z$. Thus  a braid representative
of $\delta_*$ in  the six-string braid group based on $(-3,-2,-1,1,2,3)$ and
standard generators $\sigma_1, ..., \sigma_5$
is given by $\sigma_5 \sigma_4 \sigma_5 \sigma_3 \sigma_4 \sigma_5 \sigma_2
\sigma_3 \sigma_4 \sigma_4 \sigma_1 \sigma_2 \sigma_3 \sigma_1 \sigma_2 \sigma_2$ (see figure 2.4). Thus $\delta_* =
$$A_5 A_4 A_5A_3A_4A_5A_2A_3A_4A_1A_2A_3A_1A_2A_1$.  Since
the lifts of $\delta_*$ are $\delta$ and $\tau_2 \delta$ and $C_i$
is sent to $A_i$, it follows that one of $\delta$ or $\tau_2 \delta$
is $C_5C_4C_5C_3C_4C_5C_2C_3C_4C_1C_2C_3C_1C_2C_1$. The result follows.
$\square$

\midspace{0.1cm}
\centerline{\epsfbox{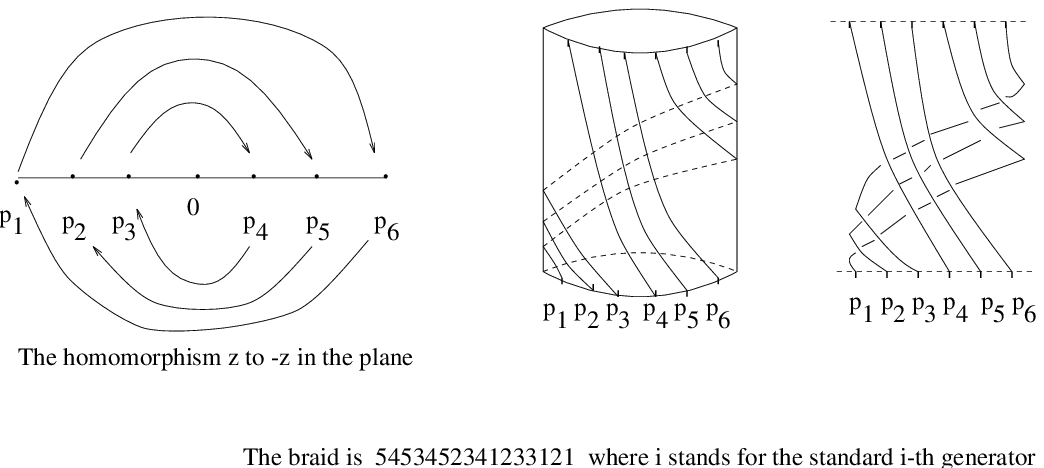}}
\midspace{0.1cm}
\centerline{Figure 2.4}

2.5. Let $S_r$ be the permutation group on $r$ boundary components
of the surface. Then there is a short exact sequence
$1 \to \Gamma_{g,r} \to \Gamma^*_{g,r} \to S_r \to1$ obtained by the
inclusion of the pure mapping class group into the mapping class group.
The lemma below shows that there are torsion elements in 
$\Gamma^*_{g,r}$ whose projections to $S_r$ generate the group.

{\bf Lemma.} \it (a) If the genus of the surface is positive,
then there are involutions in the mapping class group $\Gamma^*_{g,r}$
whose projections generate $S_r$.

(b) (see [Bi]) If the genus of the surface is zero, the mapping class group
$\Gamma^*_{0,r}$ is generated by torsion elements of orders $r-1$ and $r$.
\rm

{\bf Proof.} For part (a), note that the permutation group $S_r$ is generated
by transpositions $(ij)$. Each transposition $(ij)$ can be expressed as
a composition $\alpha \beta$ where $\alpha$ and $\beta$ are
involutions in $S_r$ having at most three fixed points. For
instance $(12) =[(12)(34) ....(n-1n)][(1)(2)(34)....(n-1n)]$
for $n$ even and $(12)=[(12)(34) ...(n-1n-2)(n)][(1)(2)(34)...(n)]$
for $n$ odd. Since any two involutions in $S_r$ with the same number
of fixed points are conjugate, it suffices to show that an involution in $S_r$
with $k \leq 3$  fixed points can be realized by an involution 
on the surface. By lemma 2.3(c), the result follows.

Part (b) follows from lemma 2.4(b).
Indeed,
the mapping class group
$\Gamma^*_{0,r}$ is generated by $A_1, ..., A_r$ (see [Bi]). 
Furthermore, any two half-twists are conjugate. Thus it suffices to
show a half-twist, say $A_{r-1}$, is a product of periodic homeomorphisms. 
By lemma 2.4(b), $A_{r-1} = (A_1...A_{r-2})^{-1}(A_1A_2....A_{r-1})$
is a product of $2\pi/(r-1)$ and $2 \pi/r$-rotations.
$\square$

2.6. 
{\bf Lemma.} \it Suppose $a$ and $b$ are two simple loops 
which intersect transversely at a point in a surface $\Sigma_{g,r}$.
Then, if the genus of the surface is one, there is a 
period-3 homeomorphism sending $a$ to $b$. If the genus of
the surface is at least two, there is a product of involutions
sending $a$ to $b$. \rm

\rm

{\bf Proof.}
Note that a regular neighborhood $N$ of $a \cup b$ is a 1-holed torus.
If the genus of the surface is at least two, we can find 
a non-separating simple loop $c$ in the surface $\Sigma_{g,r} -N$.
Now $c$ is disjoint from both $a$ and $b$ and the complements
of $a \cup c$ and $b \cup c$ in $\Sigma_{g,r}$ are connected. Thus
there are involutions in the surface sending $a$ to $c$ and $c$ to $b$
respectively. The result follows.

If the genus of the surface is one, then the completion $\Sigma^*_{1,r}$
is the torus. Furthermore, both $a$ and $b$ are still simple loops
in the torus intersecting at one point. Let $A$ and $B$ be the
Dehn-twists on these two loops. Then the 3-fold symmetry $ABAB$ sends
$b$ to $a$ in the torus. By lemma 2.3(a), there is a 3-fold symmetry
$f$ in $\Sigma_{1,r}$ inducing the 3-fold symmetry. Thus we have
$f(b) = a$.
$\square$

\S3. {\bf The Commutator Subgroup of the Mapping Class Group}

We show in this section that the commutator subgroup $[\Gamma^*_{g,r},
\Gamma^*_{g,r}]$ of the mapping class group is generated by torsion
elements in $\Gamma^*_{g,r}$.

By lemma 2.5(b), it suffices to prove the result for positive genus
surfaces.
By the short exact sequence $ 1 \to \Gamma_{g,r} \to \Gamma^*_{g,r} 
\to S_r$, proposition 2.2 and lemma 2.5(a), we conclude that
the mapping class group is generated by torsion elements and
Dehn-twists on non-separating
simple loops $a_i$'s so that either
$a_i$ is disjoint from $a_j$ or they intersect at one point.

Now if a group $G$ is generated by elements $g_i$'s, then its
commutator subgroup $[G,G]$ is normally generated by the commutators
$[g_i, g_j]$. On the other hand, if $t$ is a torsion element,
then $[s,t] = (sts^{-1})t^{-1}$ is a product of two torsions.
Thus it suffices to show that the commutator $[A,B]$
is a product of torsions in $\Gamma^*_{g,r}$ when
$a$ intersects $b$ in one point and $A= D_a$, $B = D_b$. By Dehn's relation
$ABA = BAB$, we obtain that
$[A,B]= ABA^{-1}B^{-1} = B^{-1}A$. Thus it suffices to
show that $B^{-1}A$ is a product of torsion elements.
But by lemma 2.6, there is
a product $f$ of torsion elements sending $a$ to $b$. This
implies $B = fAf^{-1}$. Thus $[A, B ] = f (A^{-1} f^{-1} A)$
is a product of torsion elements.

\S4. {\bf The First Homology of the Mapping Class Group}

We finish the proof of main theorem for surfaces of positive genus in this
section. By the short exact sequence,
$ 1 \to [\Gamma^*_{g,r}, \Gamma^*_{g,r}]
\to \Gamma^*_{g,r} \to H_1(\Gamma^*_{g,r}) \to 1$ and the result
in \S3, it suffices to show that the projections of the torsion
elements in the first homology generate the group.
To prove this, we shall recall the result on the first
homology of the pure mapping class group. Then we shall calculate
the first homology of the mapping class group and identify its
generators.

4.1. 
The first homology of the pure mapping class group $\Gamma_{g,r}$
 was shown to be torsion of order dividing 10
($g>1$) by Mumford [Mu]. Birman in [Bi1] showed that it is torsion of
order dividing 2 for $g>2$ and Powell [Po] showed the group is vanishing
for $g>2$. The most elegant proof of it is due to Harer in [Ha].

{\bf Theorem.} \it The first homology group
of the mapping class group of positive genus surface is as 
follows.
$H_1(\Sigma_{g,r}) \simeq 0, \bold Z_{10},$ and $ \bold Z_{12}$
for $g \geq 3$, $g=2$ and $g=1$ respectively. Furthermore,
a generator of the cyclic group $H_1(\Sigma_{g,r}) $ is given by
the image of the Dehn-twist on a non-separating simple loop.
\rm

A proof of it for low genus surface can be derived from the
presentation of the mapping class group obtained in [Lu].
In particular, for positive genus $g$,
 suppose $\rho$ is the group homomorphism from
$\Gamma_{g,r} \to \Gamma_{g,0}$ 
induced by
the inclusion $\Sigma_{g,r} \to \Sigma_{g,0}$. Then the induced
homomorphism $\rho_*$ between the first homology groups is an
isomorphism.

4.2. In this section, we will use the theorem 4.1 to find the
first homology group of the mapping class group $\Gamma^*_{g,r}$.

Let $\phi: \Gamma^*_{g,r} \to \Gamma_{g,0}$ be the group
homomorphism induced by the completion construction $\Sigma_{g,0}
=\Sigma^*_{g,r}$ and $\psi: \Gamma^*_{g,r} \to S_r$ be the
homomorphism induced by the action of the mapping class group
on the set of boundary components. Note that
$H_1(S_r) = \bold Z_2$ for $r \geq 2$.

{\bf Proposition.} \it 
The group homomorphism $\phi_* \oplus \psi_*: H_1(\Gamma^*_{g,r})
\to H_1(\Gamma_{g,r}) \oplus H_1(S_r)$ is an isomorphism for
positive genus surfaces. 
\rm

We will defer the proof to the last section.

4.3. We now use the proposition 4.2 to prove the remaining
part of the main theorem. Namely, for $(g, r) \neq (2, 5k+4)$,
the first homology group $H_1(\Gamma^*_{g,r})$ is generated
by the images of the torsion elements. 

First of all, by lemma 2.5 (a), there are torsion elements whose
projections in $H_1(S_r)$ generate the group. Thus it remains
to show that images of the torsion elements in $H_1(\Gamma_{g,0})$
generate the group. The result is clear for $g \geq 3$. When
$g=1$, let $A$ and $B$ be the Dehn-twists on two simple loops
intersecting at one point. Then group $H_1(\Gamma_{1,0})$ is
generated by $[A]$. Now $\tau_3 = ABAB$ and $\tau_4= ABA$ shows
that $A = \tau_3 \tau_4^{-1}$. Thus the first homology
group is generated by the images of the 3-fold symmetry $\tau_3$ and
the  4-fold symmetry $\tau_4$.
By lemma 2.3(a), both of these symmetries are induced by some symmetries
in $\Sigma_{1,r}$. Thus the result follows for $g=1$.
When the genus is two, the first homology group $H_1(\Gamma_{2,0})
\simeq \bold Z_{10}$ is generated by the image $[D]$ of a Dehn-twist on a
non-separating simple loop. Now by lemma 2.4(c), the 5-fold symmetry
$\tau_5$ and the 2-fold symmetry $\delta$ become $4[D]$ and $5[D]$
in $H_1(\Gamma_{2,0})$. Thus $H_1(\Gamma_{2,0})$ is generated by
the image of the torsion elements. By lemma 2.3(b) and (d), both
of the symmetries $\tau_5$ and $\delta$ are induced by a 5-fold
and 2-fold symmetries on $\Sigma_{2,r}$. Thus the result follows for
$g=2$. 

Finally, we need to show that the mapping class group $\Gamma^*_{2,5k+4}$
is not generated by torsions. Indeed, if it is generated by torsion
elements, then due to $H_1(\Gamma^*_{2, 5k+4}) = \bold Z_5 \oplus \bold Z_2
\oplus Z_2$, there must be torsion elements in $\Gamma^*_{2,5k+4}$
of order 5. This implies that $\bold Z_5$ acts on $\Sigma_{2,5k+4}$. 
This contradicts lemma 2.3(b).
But our proof above shows that torsion elements generate the kernel of the 
homomorphism $\theta \phi:
\Gamma^*_{2, 5k+4} \to \bold Z_5$ where $\theta: H_1(\Gamma_{2,0})
\to \bold Z_5$ is the onto homomorphism. This shows that the group
generated by torsion elements is a proper subgroup containing an
index 5 subgroup. Thus it must be the index 5 subgroup.

4.4. We prove proposition 4.2 in this section.

First of all, we note that the homomorphism $\phi_* \oplus \psi_*$
is an epimorphism. Indeed, clearly both $\phi_*$ and $\psi_*$
are epimorphisms. On the other hand, if $r \geq 2$, a half-twist
$D^{1/2}_a$ is in the kernel of $\phi_* $ and is sent to the
generator of $H_1(S_r)$ by $\psi_*$. Thus  $\phi_* \oplus \psi_*$
is onto. 

To show the kernel is trivial, we will prove that the number of
elements in  $H_1(\Gamma^*_{g,r})$ is at most
$|H_1(\Gamma_{g,0})||H_1(S_r)|$.
The group $\Gamma^*_{g,r}$ is generated by Dehn-twists on
non-separating simple loops and half-twists. Since
any two Dehn-twists on non-separating simple loops are
conjugate and any two half-twists are conjugate, the 
first homology group $H_1(\Gamma^*_{g,r})$ is generated
by a Dehn-twist $[D_a]$ on a non-separating simple loop
and a half-twist $[D^{1/2}]$. By the standard relations
on Dehn-twists [De], [Lu], we see that $[D_a]$ has order at most $n$ in
$H_1(\Gamma^*_{g,r})$ where $n = |H_(\Gamma_{g,0})|$. If
$r \leq 1$ then we see that $|H_1(\Gamma^*_{g,r})|
= |H_1(\Gamma_{g,0})||H_1(S_r)|$. If $r \geq 2$,
the composition $D^{1/2} D^{1/2}$ is  a Dehn-twist $D_b$
where $b$ is a simple loop bounding a 3-holed sphere with
two boundary components $B_1$ and $B_2$ of the surface (see figure 4.1).
We claim that the Dehn-twist $D_b$ is a product of commutator
in $\Gamma_{g,r}$. Assuming the claim, then $[D^{1/2}]$
has order 2 in the first homology. Thus there is at most
$2n$ elements in the homology group $H_1(\Gamma^*_{g,r})$.
The result follows.

\midspace{0.1cm}
\centerline{\epsfbox{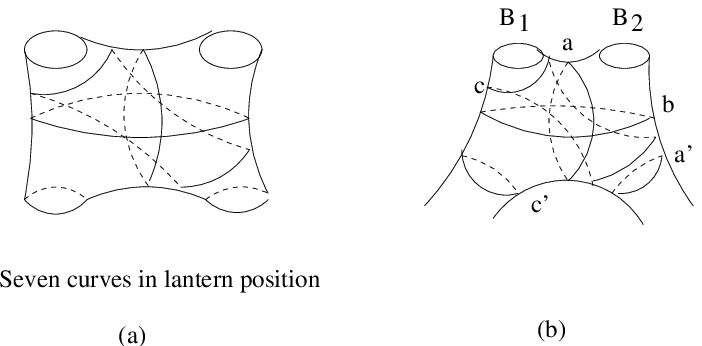}}
\midspace{0.1cm}
\centerline{Figure 4.1}

To see that $D_b$ is a product of commutators, we construct
a 4-holed sphere bounded by two non-separating simple loops
$a'$, $c'$ and $B_1$ and $B_2$. Let $a$ and $c$ be two simple
loops inside the 4-holed sphere so that $\{a,b,c\}$ forms
a lantern pattern. See figure 4.1.(b).
 Then the lantern relation
shows that
$ABC = A'C'$ where capital letters are the Dehn-twists on small
letter loops. In particular,  
$D_b = B = A'A^{-1} C'C^{-1}$. But each of $(a,a')$ and $(c,c')$
is a pair of disjoint non-separating simple loops. Thus there
are orientation preserving homeomorphisms $f$ and $h$ leaving each
boundary component invariant sending $a'$ to $a$
and $c'$ to $c$ respectively.  This shows $A' = fAf^{-1}$ and $C'
= hCh^{-1}$. Thus $D_b =[f,A][h,C]$ is a product of commutators
in the pure mapping class group.
$\square$

\bigskip

4.5. {\bf Remarks.}

1.  It can be shown that a Dehn-twist on a separating simple loop
is a product of commutators in the pure mapping class group
$\Gamma_{g,r}$ for positive genus surfaces.

2. Due to homological reason, the  period-5 and period-3 elements in the generating
set for the mapping class group of surfaces of genus 2 and 1 cannot be dropped.

3. If $\bold Z_3$ acts on a closed genus $g$ surface, 
then it has at most $g+2$ fixed
points. Indeed, if $t$ is the number of fixed points, then the 
quotient space has genus
$g'$ with $t$ many branched points of order 3. Thus the 
Euler characteristic calculation
shows:   $ 2-2g-t = 3(2- 2g' -t)$ or  the same $ t = 2+ g - 2g'$. Thus $t \leq 2+g$.
It is easy to show that there are $\bold Z_3$ actions on the genus $g$ closed surface
with $g+2$ fixed points. Here is one way to see it. Consider 
the standard $ \bold Z_3$
action on the tours with 3 fixed points. Now remove an equivariant neighborhood of
two fixed points, we obtain a $\bold Z_3$ action on the 2-holed torus with one fixed
point and  leaving each boundary component invariant. But the same argument, there
is a $\bold Z_3$ action on the  1-holed torus with 2 fixed points. Now each closed
surface can be decomposed as a union of  two 1-holed tori and $(g-2)$ 2-holed tori.
See figure 4.2.

\midspace{0.1cm}
\centerline{\epsfbox{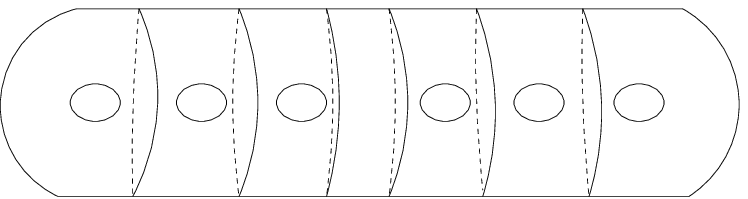}}
\midspace{0.1cm}
\centerline{Figure 4.2}

 We can glue these $\bold Z_3$ actions on the 1-holed tori and 2-holed tori to produce a $\bold Z_3$ action on the closed surface with exactly $g+2$ fixed points. By the work of Nielsen, there is only one $\bold Z_3$ action on the
closed genus $g$ surface with $g+2$ fixed points. These seems to be the analogous
action to the hyperelliptic involutions. We do not know the Dehn-twist expression
of the generator of the $\bold Z_3$ action.
Since $g+2 \geq 2$, by the same argument as in lemma 2.3, each compact surface
$\Sigma_{g,r}$ admits a $\bold Z_3$ action with at least $g$ fixed points.

\bigskip

\centerline{\bf Reference}

[Bi] Birman, J.S.:  Braids, links, and mapping class groups.
Ann. of Math. Stud., 82, Princeton Univ. Press, Princeton, NJ, 1975

[Bi1] Birman, J.S.:
Abelian quotients of the mapping class group
of a $2$-manifold. Bull. Amer. Math. Soc. 76 (1970) 147-150,
and Bull. Amer. Math. Soc. 77 (1971) 479.

[BH] Birman, J.S.; Hilden, H.:
On isotopies of
homeomorphisms of Riemann surfaces. Ann. of Math. (2) 97 (1973), 424-439.

[De] Dehn, M.: Papers on group theory and topology. J. Stillwell (eds.).
 Springer-Verlag, Berlin-New York, 1987.

[Ge] Gervais, S.: A finite presentation of the mapping class group of
an oriented surface, Topology, to appear.

[Har] Harer, J.: The second homology group of the mapping class group
of an orientable surface. Invent. Math. 72 (1083), 221-239.

[Jo] Johnson, D.:  The structure of the Torelli group. I. A
finite set of generators for $T$. Ann. of Math. (2) 118 (1983), no. 3,
42-442. 

[Li] Lickorish, R.: A representation of oriented 
combinatorial 3-manifolds. Ann.  Math.  72 (1962), 531-540

[Lu] Luo, F.: A presentation of the mapping class groups. Math.
Res. Lett. 4 (1997), no. 5, 735-739. 

[Mu] Mumford, D.:
Abelian quotients of the Teichm\"uller modular
group. J. Analyse Math. 18 1967 227-244. 

[Po]  Powell, J.: Two theorems on the mapping class group of a
surface. Proc. Amer. Math. Soc. 68 (1978), no. 3, 347-350. 

[Vi] Viro, O.: Links, two-fold branched coverings and braids,
 Soviet Math. Sbornik, 87 no 2 (1972), 216-22.

[Wa] Wajnryb, B.: 
    Mapping class group of a surface is
generated by two elements. Topology 35 (1996), no. 2, 377-383. 

\end

\end